\documentclass[12pt]{amsart}
\usepackage{epsfig,amssymb,epic,eepic,color}

\newcommand{\be}{\begin{enumerate}}
\newcommand{\ee}{\end{enumerate}}
\newcommand{\bi}{\begin{itemize}}
\newcommand{\ei}{\end{itemize}}

\def\R{\mathbb{R}}

\def\a{\alpha}
\def\om{\omega}

\def\ga{\gamma}    
\def\Ga{\Gamma}

\def\al{\alpha}
\def\be{\beta}
 
\def\de{\delta}

\def\ds{\displaystyle}

\def\nd{\noindent}
\def\bull{\hfill$\Box$}

\begin{document}
\thispagestyle{empty}%\today
\vskip 1cm
\begin{center}{\sc Essential curves in handlebodies and topological
 contractions }
\end{center}
\title{}
\author{ V. Grines}
 \address{N. Novgorod State University,
Gagarina 23, N. Novgorod, 603950 Russia.}
\email{grines@vmk.unn.ru}
\author{F. Laudenbach}
\address{Laboratoire de math\'ematiques Jean Leray, UMR 6629 du CNRS,
Facult\'e des Sciences et Techniques,
Universit\'e de Nantes,
2, rue de la Houssini\`ere,
F-44322 Nantes cedex 3, France.}
\email{francois.laudenbach@univ-nantes.fr}

\keywords{Heegaard splitting, compression disk,  
North-South diffeomorphism}
\subjclass{57M25, 37D15}
\begin{abstract}  If $X$ is a compact set, a {\it topological contraction}
 is a self-embedding 
$f$ such that the intersection of the successive images 
$f^k(X)$, $k>0$, consists of one point.
In dimension 3, we prove that 
there are smooth topological contractions of the handlebodies 
of genus $\geq 2$ whose image is essential. 
Our proof is based on an easy criterion for 
a simple curve to be essential in a handlebody. 
\end{abstract}
\maketitle
\vskip 2cm
\section{Introduction}

\medskip

\nd For a compact set $X$ and  a topological embedding $f:X\to X$,
we shall say that $f$ is a {\it topological contraction} if 
$\ds \cap_{k\geq 0}f^k(X)$ consists of one point. 
We shall show that such a 
contraction can be very complicated when $X$ is a 3-dimensional handlebody.
 Namely, we have the following  result for which some more classical
 definitions
will be recalled thereafter.

\medskip

\nd {\bf Theorem A.}
 {\it There exists a  North-South 
diffeomorphism $f$ of the 3-sphere $S^3$ and 
a Heegaard decomposition $S^3= P_-\cup P_+$ 
of genus $g\geq 2$ with the following properties:

\nd {\rm 1)} $f\vert P_+$ is a topological contraction;
 
\nd {\rm 2)}
 $f(P_+)$ is essential in $P_+$.}\\

We shall limit ourselves to $g=2$, since the generalization will be clear. 
We recall that a 3-dimensional {\it
handlebody} 
of genus 2 is diffeomorphic to the  regular
 neighborhood $P$ in $\R^3$ of the planar figure eight $\Ga$. A {\it 
compression disk} of $P$ is a smooth embedded disk in $P$ whose boundary
 lies in  $\partial P$ in which it  is not
 homotopic to a point. Among the compression disks are the 
{\it meridian} disks $\pi^{-1}(x)$, where $x$ is a regular 
point\footnote{Any point other than  the center of the figure eight.} in $\Ga$
and $\pi:P\to \Ga$ is the regular neighborhood projection (that is,
 a submersion over the smooth part of $\Ga$). A subset $X$ of $P_+$
is said to be {\it essential} in $P_+$ if it intersects every compression 
disk\footnote{ 
This definition goes back to  Rolfsen's book \cite{rolf} p. 110.}
.\\

\nd A diffeomorphism $f$ of $S^3$ is a  North-South diffeomorphism
if it  has  two fixed points only, one source $\a\in P_-$ and one sink 
$\om\in P_+$, every other orbit going from $\a$ to $\om$.\\

\nd A {\it Heegaard splitting}  of $S^3$ is made of an
 embedded surface dividing $S^3$  into two handlebodies. According to a famous 
theorem of F. Waldhausen 
such a decomposition is unique up to diffeomorphism \cite{wald1}
 (hence up to isotopy after Cerf's theorem $\pi_0(Diff_+S^3)=0$ \cite{cerf}).
 It is not hard to
 prove that the phenomenon mentioned in theorem A
does not happen with a Heegaard 
splitting of genus 1: if $T$ is a solid torus and $f$ is a topological
 contraction of $T$, then there is a compression disk of $T$
avoiding $f(T)$.\\

\nd The example which we are going to construct for proving  theorem A
is based on the next theorem, for which some more notation is introduced.
 Let $\Ga_0\subset \Ga$ be a simple closed curve. 
There exists a solid torus 
$T\subset \R^3$ which contains $P$ and which is a tubular neighborhood of 
$\Ga_0$. Let $i_0: P\to T$ be this inclusion. We say that a simple curve is 
unknotted in $T$ if it bounds an embedded disk in $T$.\\

\nd {\bf Theorem B.} {\it There exists an essential simple curve $C$ in $P$
such that $i_0(C)$ is unknotted in $T$. }\\

\nd Theorem B looks very 
easy as it is simple to draw a simple curve 
 which intuitively satisfies its  conclusion.
 Nevertheless, it appears that there are very few criteria for 
proving that a curve is essential in $P$. We are going to give one which is
not  algebraic in nature.  Question: does  there exist a topological 
algebraic tool which plays the same role.

\bigskip
\nd The second author is 
 grateful to Sylvain Gervais and 
 Nathan Habegger for interesting conversations on 
link invariants, in particular on Milnor's invariants \cite{milnor}. 
He is also grateful to the organizors\footnote{
Michel Boileau, Thomas Fiedler, John Guaschi and Claude Hayat.}
 of the conference in Toulouse in memory
of Heiner Zieschang ({\it Braids, groups and manifolds}, Sept. 2007), 
who offered him the opportunity to give a short talk on that subject.
\bigskip
\section{ Essential  curves}

\medskip

\nd Our candidate for $C$ in Theorem B is pictured in figure 1.\\ 

\begin{center}
\begin{picture}(0,0)%
\includegraphics{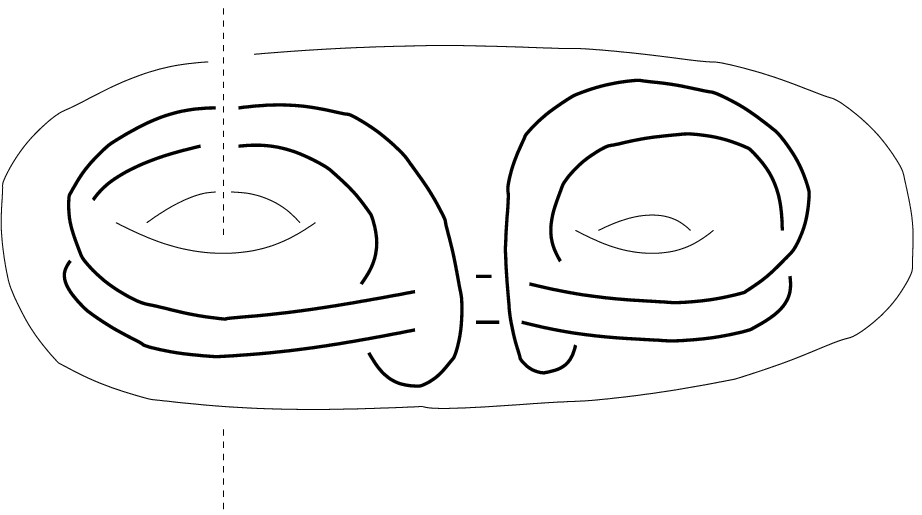}%
\end{picture}%
\setlength{\unitlength}{1934sp}%
\begingroup\makeatletter\ifx\SetFigFont\undefined%
\gdef\SetFigFont#1#2#3#4#5{%
  \reset@font\fontsize{#1}{#2pt}%
  \fontfamily{#3}\fontseries{#4}\fontshape{#5}%
  \selectfont}%
\fi\endgroup%
\begin{picture}(8956,4974)(1413,-5398)
\put(5701,-1411){\makebox(0,0)[lb]{\smash{{\SetFigFont{12}{14.4}{\familydefault}{\mddefault}{\updefault}{\color[rgb]{0,0,0}$C$}%
}}}}
\end{picture}%

Figure 1 %-- Here $\sigma$ is 1- or 2-dimensional.

\end{center}

\medskip

\nd It is clear that  $i_0(C)$ is unknotted in $T$ (or, equivalently, in the
 complement of the vertical axis which is drawn
 on figure 1 and whose 
$T$ is a compact retract by isotopy deformation). Instead of proving that $C$ 
is essential in $P$, we are going to prove a stronger result. Clearly 
Proposition 1 below implies Theorem B.\\

\nd{\bf Proposition 1.}
{\it Let $p:\widetilde P\to P$ be the universal cover of $P$ and let 
$\widetilde C$ be the preimage $p^{-1}(C)$. Then  $\widetilde C$ is essential
in $\widetilde P$.}\\

\nd {\bf Proof.} We have the following description of $\widetilde P$: it is 
a 3-ball with a Cantor set $E$ removed from its bounding
 2-sphere\footnote{Take the universal cover of $\Gamma$ properly embedded in
 the hyperbolic plane and take a 3-dimensional thickening of it.}. This Cantor
set is the set of ends of  $\widetilde P$. A simple curve in 
$\partial\widetilde P$ is not homotopic to zero if it divides $E$ into two 
non-empty parts. We get a fundamental domain $F$ for the action of $\pi_1(P)$
on $\widetilde P$ by cutting $P$ along two non-parallel meridian disks 
$D_0$ and $D_1$. Here is a description of $\widetilde C\cap F$ (see figure 2):
$F$ is a 3-ball whose boundary consists of four disks $d_0,d'_0,d_1,d'_1$
and a punctured sphere $\partial_0F$. We have $p(d_0)=p(d'_0)=D_0$
and  $p(d_1)=p(d'_1)=D_1$. We have four strands in $\widetilde C\cap F$:
$\ell_1$ and $\ell_2$ joining $d_0$ and $d_1$, $\ell'_0$ (resp. $\ell'_1$)
whose end  points belong to  $d'_0$ (resp. $d'_1$). Moreover 
$\ell'_i$, $i=0,1$, is linked with $\ell_j$, $j=1,2$, in the following sense:
any embedded surface  whose boundary is made of $\ell'_i$ and a simple 
 arc in $d'_i$ intersects $\ell_j$ for $j=1,2$ 
(the algebraic intersection number is 1 for  some choice of orientations).

\bigskip

\begin{center}
\begin{picture}(0,0)%
\includegraphics{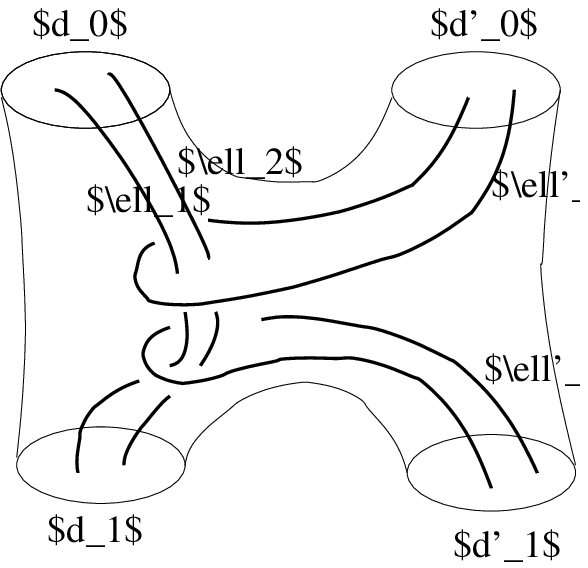}%
\end{picture}%
\setlength{\unitlength}{1934sp}%
\begingroup\makeatletter\ifx\SetFigFont\undefined%
\gdef\SetFigFont#1#2#3#4#5{%
  \reset@font\fontsize{#1}{#2pt}%
  \fontfamily{#3}\fontseries{#4}\fontshape{#5}%
  \selectfont}%
\fi\endgroup%
\begin{picture}(5649,5496)(1414,-5869)
\put(1726,-661){\makebox(0,0)[lb]{\smash{\SetFigFont{12}{14.4}{\familydefault}{\mddefault}{\updefault}{\color[rgb]{0,0,0}$d_0$}%
}}}
\put(5626,-661){\makebox(0,0)[lb]{\smash{\SetFigFont{12}{14.4}{\familydefault}{\mddefault}{\updefault}{\color[rgb]{0,0,0}$d'_0$}%
}}}
\put(1876,-5611){\makebox(0,0)[lb]{\smash{\SetFigFont{12}{14.4}{\familydefault}{\mddefault}{\updefault}{\color[rgb]{0,0,0}$d_1$}%
}}}
\put(5851,-5761){\makebox(0,0)[lb]{\smash{\SetFigFont{12}{14.4}{\familydefault}{\mddefault}{\updefault}{\color[rgb]{0,0,0}$d'_1$}%
}}}
\put(6226,-2236){\makebox(0,0)[lb]{\smash{\SetFigFont{12}{14.4}{\familydefault}{\mddefault}{\updefault}{\color[rgb]{0,0,0}$\ell'_0$}%
}}}
\put(6151,-4036){\makebox(0,0)[lb]{\smash{\SetFigFont{12}{14.4}{\rmdefault}{\mddefault}{\updefault}{\color[rgb]{0,0,0}$\ell'_1$}%
}}}
\put(2251,-2386){\makebox(0,0)[lb]{\smash{\SetFigFont{12}{14.4}{\familydefault}{\mddefault}{\updefault}{\color[rgb]{0,0,0}$\ell_1$}%
}}}
\put(3151,-2011){\makebox(0,0)[lb]{\smash{\SetFigFont{12}{14.4}{\familydefault}{\mddefault}{\updefault}{\color[rgb]{0,0,0}$\ell_2$}%
}}}
\end{picture}

Figure 2

\end{center}

\medskip
\nd Globally $\widetilde C$ looks like an infinite Borromean chain: any finite 
number of components is unlinked.
Suppose the contrary
that  $\widetilde C$ is not essential and consider 
 $\Delta$, a compression  disk
of $\widetilde P$ avoiding $\widetilde C$.
 We take it to be transversal to $\widetilde D:=p^{-1}(D_0\cup D_1)$.
Let $\mathcal C$ be the finite family of curves (arcs or closed curves)
in $\widetilde D\cap\Delta$. An element $\ga$ 
of  $\mathcal C$  is said to be
{\it innermost} if  $\ga$ divides 
$\Delta$ into two domains, one of them being a 
disk $\de$ whose interior contains 
no element of $\mathcal C$. Take such an innermost element $\ga$; 
its associated
disk $\de$ lies in $F$, up to a covering transformation, and divides $F$ into
 two 
balls $F_0$ and $F_1$.\\

\nd {\bf Lemma 1.} {\it One of the balls, say $F_0$, avoids  $\widetilde C$.
}\\

\nd{\bf Proof.} Let us consider the case when $\ga\subset d'_0$; say that 
$\ell'_0\subset F_1$. The other
 cases are very similar. Let $\al=\de\cap d'_0$. It is a simple arc dividing
$d'_0$ into two parts. Both end points of $\ell'_0$ lie in the same part since 
$\de $ avoids  $\ell'_0$. They are joined by a simple arc $\al'$ disjoint 
from $\al$. Let $\de'$ be an embedded disk bounded by $\ell'_0\cup\al'$.
This disk can be chosen disjoint from $\de$. Indeed, if $\de\cap\de'$ is not
 empty, this intersection being  transversal, by looking at an innermost 
intersection curve on $\de$ one finds an  embedded 2-sphere $S$ 
in the complement of $\widetilde C$ with one hemisphere in $\de$ and the other
 in $\de'$. As $S$ bounds a 3-ball $B_F$ in $int\,F$, which  hence is 
also disjoint from $\widetilde C$, there is an isotopy supported in a 
neighborhood of $B_F$ whose effect on $\de'$ decreases the number of 
intersection curves with $\de$. 

\nd Once $\de\cap\de'$ is empty, we have $\de'\subset F_1 $. But $\ell_1$ and 
$\ell_2$ must intersect $\de'$. Hence we have $\ell_1\cup\ell_2\subset F_1$.
Similarly,  we have $\ell'_1\subset F_1$.\bull\\

\nd  One checks easily that
 there is an isotopy of $\Delta$, supported in a neighborhood of $F_0$,
till a new compression disk having fewer intersection curves with 
$\widetilde D$ than the cardinality of $\mathcal C$. Repeating this process,
 we push
$\Delta$ into a fundamental domain, say $F$. In that position we have 
$\partial \Delta\subset \partial_0F$.
Again $\Delta$ divides $F$ into two balls and one of them,
$F_0$, avoids $\widetilde C$.
This proves that $\partial\Delta$ bounds a disk in $\partial_0 F$, namely
$F_0\cap\partial_0 F $. Hence $\Delta$ is not a compression disk.
\bull\\

\nd {\bf Remark.} We used  local linking information (namely, linking 
of strands in a fundamental 
domain of the universal covering space)
 which, as in this example, follows from usual 
linking numbers and we got a global result. This method looks very efficient. 
The general criterion is the following, where we use the same notation as 
above.\\

\nd {\bf Criterion.} {\it Let $C$ be any simple closed curve in $P$. We assume 
that there is no embedded disk $\de$ in $F$ satisfying:

\nd {\rm 1)} the boundary of $\de$ is made of two arcs $\al$ and $\be$, where
 $\al$ is an arc in $\widetilde D$ and $\be$ is an arc in $\partial\widetilde
P\cap F$;

\nd {\rm 2)} $\de$ non trivially separates the components of
 $\widetilde C\cap F$ (both components 
of $F\setminus \de$ meet $\widetilde C$).

\nd Then $C$ is essential in $P$.}

\bigskip
\section{Proof of Theorem A}

\medskip
\nd 
We recall the embedding $i_0:P\to int\,T$. We start with a curve $C$ in $P$
 which meets   the conclusion of Theorem B. We equip it 
 with its 0-normal framing
(a section in this framing is not linked with $C$ in $\R^3$) and  we choose 
an embedding $j_0:T\to P$ whose image is  a tubular neighborhood of $C$.
Let $B$ be a small ball in $int\,T$. As $C$ is unknotted in $T$, there is an
 ambient isotopy, supported in $int\,T$, deforming $i_0$ to $i_1:P\to int\,T$
such that $i_1\circ j_0(T)$ is a standard small solid torus in $B$.
One half of the desired Heegaard splitting of genus 2 will be given by $
P_+:=i_1(P)$. At the present time $f$ is only defined on $T$ by $f:=
i_1\circ j_0:T\to int\,T$. If we compose $i_1$ with  a sufficiently 
strong contraction of $B$ into itself, then $f$ is a contraction in the
 metric sense. Hence $\ds\cap_{k>0}f^k(T)$ consists of one point.\\

\nd Choose a round ball $B'$ containing $T$ in its interior. Since
$f|T$ is isotopic to the inclusion $T\hookrightarrow \R^3$, $f$
extends as a diffeomorphism $B'\to B$, and further as a diffeomorphism
$S^3\to S^3$. We are free to choose $f:S^3\setminus B'\to S^3\setminus B$
as we like. If we compose $f^{-1}$ with a strong contraction
of $S^3\setminus B'$, the intersection  $\ds\cap_kf^{-k}(S^3\setminus B')$   
consists of one point. We now have a North-South diffeomorphism $f$ of $S^3$
which induces a topological contraction of $T$. Since $f(T)\subset int\,P_+
\subset P_+\subset int\,T$, $f$ also induces a topological contraction 
of $P_+$.\\

\nd It remains to prove that $f(P_+)$ is essential in $P_+$. We know that 
$i_1(C) $ is essential in $P_+$. As a consequence, any compression disk 
$\Delta$ of $P_+$ crosses $f(T)$. We can take $\Delta$ to be  transversal to 
$f(\partial T)$ such that no intersection curve is null-homotopic in 
$f(\partial T)$. Let $\ga$ be an intersection curve which is {\it innermost}
in  $\Delta$  and let $\de$  be the disk that $\ga$ bounds in $\Delta$.\\

\nd {\bf Lemma 2.} {\it We have $\de\subset f(T)$.}\\

\nd {\bf Proof.} If not, we have $\de\subset P_+\setminus f(int\,T)$ and the
 simple curve $\ga$ in $f(\partial T)$ is unlinked with the core 
$i_1(C)$. Therefore, up to isotopy in  $f(\partial T)$, it is a section of
 the 0-framing. In that case, $i_1(C)$ itself bounds an embedded disk in $P_+$.
This is impossible, as $i_1(C)$ is essential in $P_+$.\bull\\

\nd Therefore $\de$ is a compression disk of the solid torus $f(T)$.
But $P_+=i_1(P)$, like $P$ itself, is essential in $T$. Hence $f(P_+)$
is essential in $f(T)$ and $\de$ must cross $f(P_+)$. \bull

${}$\\

\end{document}